\documentclass[11pt]{amsart}
\usepackage[T1]{fontenc}
\usepackage[latin9]{inputenc}
\usepackage{amssymb}
\usepackage{amsmath}
\usepackage{graphicx}
\usepackage{enumerate}
\usepackage{framed}

\usepackage[english]{babel}

\newcommand{\MSOL}{\mathrm{MSOL}}
\newcommand{\SOL}{\mathrm{SOL}}
\newcommand{\HOL}{\mathrm{HOL}}
\newcommand{\FPT}{\mathrm{FPT}}

\newcommand{\Z}{\mathbb{Z}}
\newcommand{\R}{\mathbb{R}}
\newcommand{\C}{\mathbb{C}}

\newcommand{\bP}{\mathbf{P}}
\newcommand{\bNP}{\mathbf{NP}}
\newcommand{\bPH}{\mathbf{PH}}
\newcommand{\bPSpace}{\mathbf{PSpace}}
\newcommand{\bVP}{\mathbf{VP}}
\newcommand{\bVNP}{\mathbf{VNP}}

\newtheorem{thm}{Theorem}

\newif\ifskip
\skiptrue

\begin{document}
\title{How I got to like graph polynomials}
\author{J.A. Makowsky}
\address{Faculty of Computer Science, \newline Technion--Israel Institute of Technology, Haifa, Israel}
\email{janos@cs.technion.ac.il}
\thanks{I would like to thank J. Kirby for various suggestion on how to improve the paper. }



\date{August, 05, 2023}

\dedicatory{For Boris Zilber on his 75th birthday}

\begin{abstract}
I trace the roots of my collaboration with Boris Zilber, which combines categoricity theory, finite model theory,
algorithmics, and combinatorics.
\end{abstract}
\maketitle

\section{Introduction and dedication}
Boris Zilber played a crucial r\^{o}le in my work on graph polynomials.
Some of my work he inspired and in which he contributed, is summarized in  Kotek et al.\ \cite{ar:KotekMakowskyZilber11}.
A preprint was posted as Makowsky and Zilber \cite{makowsky2006polynomial}
and a conference paper was published as Kotek et al.\ \cite{kotek2008counting}.
These are our only published joint papers. 
Since then a general framework for studying graph polynomials
has slowly evolved. 
It bears witness to the impact of Boris on my own work.
In this paper I will sketch how I got to like graph polynomials.
Boris and I both started our scientific career in model theory. 
Boris pursued his highly influential work in various directions of infinte model theory.
My path towards graph polynomials 
took a detour into the foundations of computer science, only to lead me back to model theoretic methods
in finite combinatorics. I sketch here how, step by step, I ended up discussing graph polynomials with Boris. 
Some of those steps owe a lot to serendipity, others were triggered by natural questions arising from
previous steps. These steps are described in Sections \ref{se:2}-\ref{se:6}.
Sections \ref{se:7}-\ref{se:8} describe some of the original ideas underlying the model-theoretic approach to graph polynomials.
Section \ref{se:9} summarizes where this encounter with Boris has led me.
Ultimately, it looks as if Boris' influence on my path was inevitable, but only in retrospect.
Meeting Boris in Oxford was a chance encounter with unexpected consequences.
I would like to thanks Boris for giving me an important and fruitful impulse.
Happy birthday, and many years of productive mathematics to come, till 120.
\section{Morley's 1965 paper}
\label{se:2}
My first attempt to tackle open problems in model theory was a consequence of reading 
M. Morley's
fundamental paper 
\cite{ar:Morley65} 
on categoricity in power,
in the {\em undergraduate seminar}
in mathematical logic at ETH Z\"urich in 1969. 
The seminar was
held by E. Specker and H. L\"auchli, 
and regularly attended by the then still very lucid octogenarian P. Bernays.

Building on  earlier work by A. Mostowski, A. Ehrenfeucht and R. Vaught,
M. Morley proved in 1965
a {\bf truly deep } theorem in model theory:
\begin{thm}[Morley's Theorem]
Let $T$ be a first order theory and assume $T$ has no finite models and  
is $\kappa$-categorical for {\bf some} uncountable $\kappa$. 
Then $T$ is $\kappa$-categorical for {\bf every} uncountable $\kappa$.
\end{thm}
More importantly, even, the paper ended with a {\bf list of questions} and
many logicians and mathematicians were attracted by these. Among them I remember
J. Baldwin and A. Lachlan, J.-P. Ressaire, D. Lascar, M. Makkai, V. Harnik, S. Shelah, B. Zilber, M.A. Taitslin and
his school, see Taitslin \cite{Taitslin1970},  and myself. 
In my MSc thesis from 1971 \cite{makowsky1974some}, I managed to prove the following:
\begin{thm}
\begin{enumerate}[(i)]
\ 
\item
A first order theory $T$ which is $\aleph_0$-categorical and strongly minimal (hence
categorical in all infinite $\kappa$) cannot be finitely axiomatizable.
\item
There is a finitely axiomatizable complete first order theory $T$ without finite models which is superstable.
\item
If there is an 
{\em infinite, finitely presentable group  $G$ with only finitely many conjugacy classes},
there is also a complete finitely axiomatizable $\aleph_1$-categorical theory $T_G$ without finite models.
\end{enumerate}
\end{thm}
After I finished my MSc thesis, E. Specker drew my attention to a Soviet paper by 
Mustafin and Taimanov \cite{Mustafin1970countable}, and as a result of this,
I started a correspondence with its second author.
Before 1985 there were very few authors citing my work.
Among them G. Ahlbrandt, a PhD-student of J. Baldwin,
P. Rothmaler, P. Tuschik from the German Democratic Republic,
and
B. Zilber, M. Peretyat'kin, and A. Slissenko from the Soviet Union.
Boris was one of the first to notice and cite my work on categoricity.
I soon realized that I could not make any further progress on these questions.
I had no new ideas, and competition was overwhelming. S. Shelah's sequence of papers inspired by these
open questions led many young researchers to abandon this direction of research in model theory.
The finite axiomatizability questions were finally solved
by Peretyat'kin \cite{peretyatkin1980example} and Zilber \cite{zilber1981totally}. 
M. Peretyat'kin constructed a finitely axiomatizable theory
categorical in $\aleph_1$ but not in $\aleph_0$. B. Zilber showed that no finitely axiomatizable
totally categorical first order theory exists. 
An alternative proof of this was given by G. Cherlin et al.\ \cite{cherlin1985omega}.

My first acquaintance with Boris Zilber happened via the literature.
But our paths diverged (not in the yellow wood), and we did meet personally, but not very often.

\section{From abstract model theory to computer science and graph algorithms}
\label{se:3}
After leaving Morley-type model theory,
I first worked in abstract model theory, and then in theoretical computer science. 
In computer science
I dealt with the foundations of database theory and logic programming, which led me to finite model theory.
My main tools from model theory were pebble games and the Feferman-Vaught theorem and its generalizations.
Around this time I met B. Courcelle and became aware of the Robertson-Seymour theorems
and their applications to graph algorithms described by Fellows \cite{fellows1989robertson}. But it was Courcelle
\cite{courcelle1992monadic}
who first observed that logical methods would give even more applications, 
Courcelle's work on the monadic second order theory of graphs is summarized in the monumental
monograph from 2012  by Courcelle and Engelfriet \cite{courcelle2012graph}.

Let $d(G)$ be a graph parameter and $P$ be a graph property.
If deciding whether a graph $G$ on $n$ vertices with $d(G)=t$ is in $P$ can be done in time $c(t)\cdot n^s$ for some fixed $s$
which does not depend  on $d(G)$, nor on the number of vertices of $G$,
we say that {\em $P$ is Fixed Parameter Tractable ($\FPT$)}.
This concept was introduced by Fellows and Downey \cite{downey2013fundamentals}.

\begin{thm}[Courcelle, 1992]
Let $C$ be a class of finite graphs of tree-width at most $t$
and let $P$ be a graph property definable in Monadic Second Order Logic $\MSOL$.
Then checking whether a graph $G ,\in C$ with $n$ vertices is in $P$ is in $\FPT$,
in fact, it can be solved in linear time  $c(t)n$.
\end{thm}

In the mid-1990s two students were about to change my research dramatically.
My former MSc {\em Udi Rotics} returned from his experience in industry.
His MSc thesis dealt with the logical foundation of databases.
However, now he wanted to work for a PhD in {\em Graph Algorithms} but {\em without involving Logic}.
He proposed to extend the notion of {\em tree-width} of a class of graphs as a graph parameter in order
to get a new width parameter which one can
use for {\em Fixed Parameter Tractability}. 
Finally, {\em but still using Logic ($\MSOL$)},  we came up with a notion roughly equivalent to {\em clique-width},
introduced recently by Courcelle and  Olariu \cite{courcelle1994clique}.
This led to my intensive collaboration with Courcelle and Rotics
\cite{courcelle1998linear,courcelle2000linear,courcelle2001fixed}.
In my own paper \cite{makowsky2004algorithmic} I examine the algorithmic uses of the Feferman-Vaught theorem
for Fixed Parameter Tractability.
Applications of my work with B. Courcelle and U. Rotics are well summarized in Downey and Fellows \cite{downey2013fundamentals}.

In 1996 I started to supervise an immigrant student from the former USSR, {\em Gregory Kogan},
who wanted to work on the complexity of computing the {\em permanent}.
He came with impressive letters of recommendation.
He had some spectacular partial results for computing permanents of matrices over a field of characteristic $3$.
He was a virtuoso in combinatorial linear algebra.
Unfortunately, he dropped out before finishing his PhD.
M. Kaminski and I wrote up his results, published under his name alone as Kogan \cite{DBLP:conf/focs/Kogan96}.

\section{Computing permanents}
\label{se:4}
I first came across the problem of computing the permanent
at Specker's 60th birthday conference in 1980.
The permanent of an $(n \times n)$-matrix $A=(A_{i,j})$ is given as
$$
per(A) = \sum_{s:[n] \rightarrow [n]} \prod_{i \in [n]}  A_{i,s(i)}
$$
where $s$ ranges over all permutations of $[n]$.

The complexity class $\sharp\mathbf{P}$ is the polynomial time counting class.

The class of $\sharp\mathbf{P}$ consists of function problems of the form ``compute $f(x)$'', 
where $f$ is the number of accepting paths of a nondeterministic Turing machine running in polynomial time. 
Unlike most well-known complexity classes, it is not a class of decision problems but a class of function problems. 
The most difficult representative problems of this class are $\sharp\mathbf{P}$-complete. 
Counting the number of satisfying assignments for a formula of propositional logic
is $\sharp\mathbf{P}$-complete.
Typical examples would be described as follows: Let $k$ be a fixed integer. Given 
an input graph $G$ on $n$ vertices,
compute the number of proper $k$-colorings of $G$.
For $k=1,2$ this can be computed in polynomial time, but for $k \geq 3$, this is $\sharp\mathbf{P}$-complete
with respect to $\bP$-time reductions. In general $\sharp\mathbf{P}$ lies between the polynomial hierarchy $\bPH$
and $\bPSpace$, see Papadimitriou \cite{bk:papadimitriou94}.

Valiant's complexity classes $\bVP$ and $\bVNP$ are the analogues of $\bP$ and $\bNP$
in Valiant's model of algebraic computation. B\"urgisser's book \cite{burgisser2000completeness} is entirely dedicated
to this model of computation.
It is still open whether $\bP=\bNP$, and also whether $\bVP=\bVNP$.

Valiant presented the complexity classes $\mathbf{VP}$ and $\mathbf{VNP}$ at 
Specker's 60th birthday conference.

\begin{thm}[L. Valiant \cite{valiant1979complexity}]
Computing the permanent of a $\{0,1\}$-matrix is hard in the following sense:
\begin{enumerate}[(i)]
\item It is $\sharp\mathbf{P}$-complete in the Turing model of computation, and
\item $\mathbf{VNP}$-complete in Valiant's algebraic model of computation.
\end{enumerate}
\end{thm}
Valiant \cite{ar:Valiant80}
published in the proceedings of Specker's 60th birthday conference in 1980.
B\"urgisser's monograph \cite{burgisser2000completeness} 
explores Valiant's approach to algebraic complexity further.

G. Kogan studied the complexity of computing the permanent over fields of characteristic $3$
for matrices $M$ with  $rank(MM^{tr} - \mathbf{1}) = a$. He showed that for $a \leq 1$ this is easy,
and for $a \geq 2$ this is hard.

I wanted to use results from Courcelle et al.\ \cite{courcelle2001fixed} to prove something about permanents
G. Kogan could not prove.
I looked at adjacency  matrices of graphs of fixed tree-width $t$. 
Barvinok \cite{barv:96}
also studied the complexity of computing the permanent for special matrices.
He looked at matrices of fixed rank $r$.
Our results were:
\begin{thm}[A. Barvinok, 1996]
Let $\mathcal{M}_r$ be the set of real matrices of fixed rank $r$.
There is a polynomial time algorithm $\mathcal{A}_r$
which computes $per(A)$
for every $A \in \mathcal{M}_r$ 
\end{thm}
\begin{thm}[JAM, 1996]
\label{jam-per}
Let $\mathcal{T}_{w}$ be the set of adjacency matrices of graphs of tree-width at most $w$.
There is a polynomial time algorithm $\mathcal{A}_w$
which computes $per(A)$
for every $A \in \mathcal{T}_w$.
\end{thm}
The two theorems are incomparable.
There are matrices of tree-width $t$ and arbitrary large rank, and
there are matrices of rank $r$ and arbitrary large tree-width.

However, I realized that the proof of my theorem  had nothing to do with permanents.
It was much more general and really worked quite generally. It only depended on some 
{\em logical restrictions} for polynomials
in indeterminates given by the entries of the matrix.
If the matrix $A=A_G$ is the adjacency matrix of a graph $G$ where the non-zero entries are $x$, 
the permanent $per(A_G)$ can be viewed as a graph polynomial in the indeterminate $x$.
Alas, at that time I had no clue how to find  many interesting examples.

\section{From knot polynomials to graph polynomials}
\label{se:5}
During a sabbatical at ETH in Zurich I met {\em V. Turaev}, who, among other things, is an expert in knot theory. 
I showed him my Theorem \ref{jam-per}.
He suggested I should try to prove the same for the {\em Jones polynomial} from {\em Knot Theory}.
So I studied Knot Theory intensively for a few months. While visiting the Fields Institute in 1999
I attended a lecture by J. Mighton\footnote{
John Mighton is a Canadian mathematician, author, and playwright.
\\
https://en.wikipedia.org/wiki/John\_Mighton
} who lectured about the Jones polynomial for
series-parallel knot diagrams, cf. his PhD thesis, Mighton \cite{mighton2000knot}.
He showed that in this case the Jones polynomials is computable in polynomial time.
Series-parallel graphs are exactly the graphs of tree-width $2$. It seemed reasonable that the same would hold
for graphs of tree-width $k$.
Indeed, after quite an effort I proved in 
\cite{makowsky2001coloured,makowsky2005coloured}:
\begin{thm}
Assume $K$ is a knot with knot diagram $D_k$ of tree-width $k$.
Then evaluating the Jones polynomial $J(D_k;a,b)$ for fixed complex numbers $a,b \in \C$ and $D_k$ with $n$ vertices
is in $\mathbf{FPT}$ with parameter $k$.
\end{thm}
Jaeger et al.\ \cite{jaeger1990computational}a showed that
Without the assumption on tree-width evaluating the Jones polynomial is 
$\sharp\mathbf{P}$-complete for almost all $a,b \in \C$. 
Lotz and Makowsky \cite{lotz2004algebraic} analyze the complexity of the Jones polynomial in
Valiant's model of computation.

However, again the proof seemed to work for other graph polynomials as well,
among them the {\em Tutte polynomial, chromatic polynomial, characteristic polynomial, matching polynomial}.
Univariate graph polynomials are graph invariants which take values in a polynomial ring, usually 
$\Z[X]$, $\R[X]$ or $\mathbb{C}[X]$.
The univariate chromatic polynomial $\chi(G;k)$ of a graph counts the number of proper colorings of a graph
with at most $k$ colors. It was introduced by G. Birkhoff in 1912 in a  unsuccessful attempt to prove the four colour conjecture.
The characteristic polynomial of a graph is the characteristic polynomial (in the sense of linear algebra)
of the adjacency matrix $A_G$ of the graph $G$, see the monographs by Chung \cite{chung1997spectral} and
by Brower and Haemer \cite{brouwer2011spectra}.
The coefficients of $X^k$ of the matching polynomial count the number of $k$-matchings of a graph $G$,
see Lov\'asz and Plummer \cite{lovasz2009matching}. 
Both, the characteristic and the matching polynomial, have found applications in 
{\em theoretical chemistry}  as described by Trinajstic  \cite{trinajstic2018chemical}.
There are also multivariate graph polynomials.
The Tutte polynomial is a bivariate generalization of the chromatic polynomial. Both of them are widely studied, see
Dong et al.\ \cite{dong2005chromatic} and  the handbook of the Tutte polynomial edited by
Ellis-Monaghan and Moffatt \cite{ellis2022handbook}.
Other widely studied graph polynomials are listed in Makowsky \cite{makowsky2008zoo}.
However, I had no idea,
how to find {\em infinitely many natural and interesting examples}?

\section{Boris, deus ex machina}
\label{se:6}
In 2005, while attending CSL, the European Conference in Computer Science Logic in Oxford, I paid a visit to Boris Zilber,
whom I knew and had met before due to our work on Morley's problem on finite axiomatizability
of totally categorical theories.
After a few friendly exchanges the following dialogue evolved:
\begin{description}
\item[Boris] What do you work on nowadays?
\item[Me] Graph polynomials.
\item[Boris] What polynomials?
\end{description}
It seemed Boris had never heard of graph polynomials.
I gave him the standard examples (Tutte, chromatic, matching).
He immediately saw them as examples which are interpretable in some
totally categorical theory. 
I could not believe it.

We spent the next days together, verifying that all the known graph polynomials fit into Zilber's framework.
It was indeed the case. We also produced generalizations of chromatic polynomials, some of which
I later called  {\em Harary polynomials} \cite{herscovici2020harary, herscovici2020harary-a}.
They are generalizations of the chromatic polynomial based on conditional colorings introduced in
Harary \cite{pr:Harary85} in 1985. Conditional colorings are defined using
a graph property  $P$.
A $P$-coloring  $f$ of $G$ with at most $k$ colors is a function
$f: V \rightarrow [k]$ such that for every color $j \in [k]$ the set 
$ f^{-1}(j)$ induces a graph in $P$. 
Conditional colorings are studied in the literature, e.g., by Brown and Corneil \cite{brown1987generalized}, 
mostly in the context of extremal graph theory.
However, nobody wrote about the fact that
counting the number of such colorings with at most $k$ colors defines a polynomial in $k$.
The so called Harary polynomial $\chi_P(G;k)$ counts the number of $P$-colorings of $G$
with at most $k$ colors.

\section{Why is the chromatic polynomial of a graph a polynomial?}
\label{se:7}
Let $G=(V,E)$ be a graph. A proper coloring of $G$ with at most $k$ colors is a function
$f: V \rightarrow [k]$ such $f(v) =f(v')$ implies that $\neg E(v,v')$.
We think of $[k]$ as a set of colors.
In other words, if two vertices have the same color they are not adjacent.
We denote by $\chi(G;k)$ the number of proper colorings of $G$ with at most $k$ colors.
Birkhoff's proof that $\chi(G;k)$ is a polynomial in $\mathbb{Z}[k]$  uses 
deletion and contraction of edges.
Let $e =(u,v)$ be an edge of $G$. 
$G_{-e}$ is the graph $G_{-e}= (V, E - \{(u,v)\})$ where $e$ is deleted from $E$.
$G_{/e}$ is the graph $G_{/e}= (V_{/e}, E|_{V_{/e}})$ where $e$ is contracted to a single vertex to form
$V_{/e}$ and $e$ is omitted from $E$.
$f$ is a proper coloring of $G_{-e}$ if either it is a proper coloring of $G$ and $f(u) \neq f(v)$
or it is a proper coloring of $G_{/e}$ and $f(u) = f(v)$.
Furthermore,  $\chi(G;k)$ is multiplicative, i.e., if $G$ is the disjoint union of $G_1$ and $G_2$,
then 
$$\chi(G_1 \sqcup G_2;k) = \chi(G_1;k) \cdot\chi(G_2;k).$$ 
Let 
$E_n$ be
the edgeless graph  with $n$ vertices and
$E= \emptyset$. We have $\chi(E_n;k) =k^n$, and 
$$
\chi(G_{-e};k) = \chi(G;k) + \chi(G_{/e};k).
$$
By showing that one can compute $\chi(G;k)$
by successively removing edges, and this is independent of the order of the edges,
one concludes that $\chi(G;k)$ is a polynomial in $k$.
The disadvantage of this elegant proof is, that it does not generalize.

Another way of proving that $\chi(G;k)$ is a polynomial in $k$ is by noting
that for graphs on $n$ vertices we have
$$
\chi(G;k) = \sum_{i=1}^n c_i(G) k_{(i)}
$$
where 
the coefficient $c_i(G)$ is the number of proper colorings of $G$ with exactly $i$ colors and
$$
k_{(i)}= k \cdot (k-1) \cdot \ldots \cdot (k-i+1)  = \prod_{i=0}^i (k-i) = {k \choose i}\cdot i!
$$ 
the falling factorial. Note that ${k \choose i} = 0$ for $i > k$.
As $k_{(i)}$ is a polynomial in $k$ and  $\chi(G; k)$ is a sum of $n$ polynomials in $k$, the result follows.
However, this proof {\em does generalize}, and it works for all Harary polynomials.

Later  I discussed 
Zilber's view of graph polynomials with A. Blass.
We noted that for most of the graph invariants from the literature,
proving that they were polynomial invariants via
totally categorical theories was an overkill.
This led  me to formulate a 
considerably simplified approach, which indeed covered all the know examples
of graph polynomials in the literature.
This approach is a simplification of Boris' proof. It generalizes also to other types of graph polynomials such as
the bivariate Tutte polynomial and the trivariate edge elimination polynomial from 
Averbouch et al.\ \cite{averbouch2010extension}.
More intrinsic examples are also discussed in Makowsky and Zilber \cite{makowsky2006polynomial} and
Kotek et al.\ \cite{ar:KotekMakowskyZilber11}.
However, the polynomial graph invariants hidden in totally categorical theories are the most general graph invariants
which are definable in Second Order Logic $\SOL$, and even in Higher Order Logic $\HOL$, 
over the graph $G$, see  Makowsky and Zilber 
\cite[Corollary C and Theorem 3.15]{makowsky2006polynomial}. 
Furthermore, it applies to $\HOL$-definable polynomial invariants over arbitrary finite first order structures 
for finite vocabularies, rather than just to graphs. 

\section{The model-theoretic approach to the chromatic polynomial}
\label{se:8}
The way Boris looked at the chromatic polynomial was even more general.
Given a graph $G$, Boris had in mind an infinite first order structure $\mathfrak{M}(G)$ 
with universe $M$,
and a formula $\phi(x)$
such that
\begin{enumerate}[(i)]
\item
The first order theory $T(\mathfrak{M}(G))$  of $\mathfrak{M}(G)$ is totally categorical and strongly minimal
with a strongly minimal infinite set $X$ of indiscernibles.
\item
$T(\mathfrak{M}(G))$ has the finite model property, i.e., the algebraic closure $\mathrm{acl}(Y)$  in $\mathfrak{M}(G)$
of a finite subset $Y \subset X$ is finite.
\item
$\mathfrak{M}(G) \models \phi(x)$ iff $x$ is a proper coloring of $G$.
\end{enumerate}
Such theories were at the heart of his work \cite{Zilber-uct}.
From Zilber's analysis of totally categorical theories  we get in the spirit of
\cite[Theorem 1.5.5]{Zilber-uct}:
\begin{thm}
Let $\mathfrak{M}(G)$ and $X$ as above.
For every finite set $Y \subset X$ of cardinality $k$ the cardinality of the the set
$$
\{ x \in \mathrm{acl}(Y) : \mathfrak{M}(G) \models \phi(x) \}
$$ 
is a polynomial in $k$.
\end{thm}

In the case of the chromatic polynomial this looks as follows:
  
\begin{enumerate}[(i)]
\item
Let $G = (V,E)$ be a finite graph, with $|V| = n$.
\item 
Let $\mathfrak{M}(G) =(V,X;E, \bar{v})$ be a 2-sorted language with sorts $V, X$, 
a binary relation $E$ on $V$ for the edge relation, 
and $n$ constant symbols $v_1,\ldots,v_n$ of sort $V$.
\item 
The describing axioms state that $(V,E)$ is exactly the finite graph we started with, and the vertices are exactly the $v_i$.
Then a model of the axioms is specified up to isomorphism by the cardinality of $X$. 
\item
If we add axioms stating that $X$ is infinite then the first order theory $T(\mathfrak{M}$ is totally categorical. 
\item
We also get finite models $M_k$ where $|X| = k$ for each natural number $k$, 
and they are algebraically closed subsets of the infinite model $\mathfrak{M}$.
\item
We regard $X$ as a set of colors, and we identify $X^n$ with the colorings of vertices, that is, 
the set of functions $V \to X$, by identifying $f = (x_1,\ldots,x_n) \in X^n$ with the function $f(v_i) = x_i$.
\item
$f$ is a \emph{proper vertex coloring} if any two adjacent vertices have different colors. 
So the set of proper vertex colorings is defined as a subset of $X^n$ by the formula $\phi(\bar{x})$ given by
$$
\phi(\bar{x}): 
\bigwedge_{\{(i,j) : E(v_i,v_j)\}} x_i \neq x_j.
$$
\item 
The chromatic polynomial for the graph $G$ is then $\chi(G;k) = |\phi(M_k)|$. 
\end{enumerate}

Boris also showed me at our first encounter how the bivariate matching polynomial and the Tutte polynomial
can be cast in this framework. For the characteristic polynomial $p(G;x)$ the situation is a bit
more complicated, because its original definition uses the characteristic polynomial of the
adjacency matrix $A(G)$ of $G$. However, there exists a purely graph theoretic description
of the coefficients of $p(G;x)$ by Godsil \cite{bk:Godsil93}, which allows to cast $p(G;x)$ into this framework.
We note that it may be unexpectedly tricky to put a graph invariant into Boris' framework, even if one
already knows that it is a polynomial invariant.

The most general version of this can be found in Kotek et al.\ \cite[Section 8]{ar:KotekMakowskyZilber11}.
Using this method, any multivariate polynomial graph invariant deinable in $\HOL$ can be captured in this way.
In the last ten years J. Ne{\v{s}}et{\v{r}}il and his various collaborators
(A. Goodall, D. Garijo and P. Ossona de Mendez) were exploring various ways
to define such graph invariants. However, they did not reach the same generality, cf.
\cite{garijo2011distinguishing,goodall2016strongly,garijo2016polynomial}. 
The potential of the general approach as described in \cite[Section 8]{ar:KotekMakowskyZilber11}
still has not been explored in depth. It seems that its abstract generality makes it difficult
for the combinatorics community to see through this construction. On the other side,
model theorists seemingly are not interested in combinatorial applications of model theory.
Exceptions may be in extremal combinatorics, as initiated by Razborov \cite{razborov2007flag,razborov2013flag} 
and  surveyed by Coregliano and Razborov \cite{coregliano2020semantic}.
Another direction is counting the number $S_P(n)$ of graphs on $n$ vertices in a hereditary graph property $P$, 
initiated by Scheinerman \cite{scheinerman1994size}
and further pursued by Balogh et al.\ \cite{ar:BaloghBollobasWeinreich2000}
and Laskowski and Terry \cite{laskowski2022jumps}.

\section{Towards a general theory of graph polynomials}
\label{se:9}
For the last 20 years I was studying graph polynomials
\cite{makowsky2006zoo,makowsky2008zoo}, aiming to understand what they have in common.
\begin{itemize}
\item
Graph polynomials can be studied to obtain {\em information on graphs}.
\\
As an example: Evaluations of the  Tutte polynomial encodes many graph invariants.
\item
Graph polynomials can be studied as {\em polynomials indexed by graphs}.
\\
As an example: The acyclic matching polynomial of paths, cycles, 
complete graphs, and complete bipartite graphs are the
Chebyshev polynomials of the second and first kinds, Hermite polynomials, and
Laguerre polynomials, respectively.
\item
Two graph polynomials have the same {\em distinctive power} if they distinguish between the same graphs.
\end{itemize}

With my various collaborators I managed to create a new field in graph theory with
two Dagstuhl Seminars (16241, 19401), two MATRIX Institute programs, two special sessions at AMS meetings, and one  
SIAM mini-symposium.


Without Boris Zilber's eye opener I would not have pursued this line of research as far as I did.
\begin{center}
\large
Thank you, Boris!
\end{center}


\end{document}